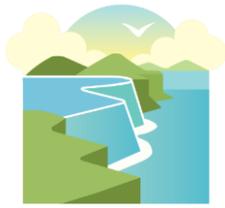
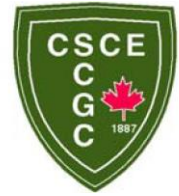

CSCE 2021 Annual Conference
*Inspired by Nature – Inspiré par la Nature*

26-29 May 2021

# MODELLING OF COUPLED SURFACE AND SUBSURFACE WATER FLOWS

Karjoun, H.[1, 3], and Beljadid, A.[1, 2, 4]
[1] Mohammed VI Polytechnic University, Green City, Morocco.
[2] University of Ottawa, Ottawa, Canada.
[3] hassan.karjoun@um6p.ma
[4] abeljadi@uottawa.ca

**Abstract:** In this paper, we propose to use the HLL finite volume scheme combined with implicit techniques for modelling the coupled surface and subsurface water flows. In our approach, we used the shallow water equations modelling surface water flow with different source terms such as variable bottom topography, friction effect, precipitation and infiltration. For subsurface water flow, the Green-Amp equation is used to simulate the infiltration process through soils. For solving the resulting nonlinear-coupled system of shallow water flow and the Green-Ampt infiltration equations, the HLL finite volume schemes with linear reconstructions of the solutions at the discrete level are implemented in order to achieve the second-order accuracy of the scheme. Appropriate discretization techniques are used for the source terms to guarantee the well-balanced property of our numerical scheme. Numerical experiments are performed to test the capability of the developed numerical scheme to simulate the coupled surface and subsurface water flows.

## 1   INTRODUCTION

Modelling of coupled systems of surface and subsurface water flows has many applications in hydrology, water management and agricultural engineering (Taccone, et al. 2020, Fernández-Pato, Gracia and García-Navarro 2018, Tatard, et al. 2008). These coupled models allow understanding the dynamic processes of surface water flows and their interactions with the infiltration process in soils, and can be used for mass balance estimation. The water mass balance analysis is used to quantify the amount of surface water for water management purposes (Pham, et al. 2016). In our study, the physically based model used for surface water flows is the shallow water equations (SWEs), also called the Saint Venant system (De St Venant 1871). The SWEs are widely used to study systems of shallow water flows where the water depth is much smaller than the horizontal length scale of these systems and the variations of the flow in the vertical direction is negligible compared to the variations in the horizontal one (Mugler, et al. 2011, Cea and Vazquez-Cendon 2012, Yoon and Kang 2004). For subsurface water flow, the Green-Ampt model (Green and Ampt 1911) is considered to describe the infiltration process through soils. This simplified physically based hydrodynamic model can be used to predict infiltration capacity, ponding time, and cumulative infiltration for layered soils.

Different physical processes occur during rainfall event including the infiltration process in soils, the dynamics of flows over variable topography with friction effect and their interactions, which make their simulations very challenging and require robust numerical techniques. In the present study, to better capture these physical processes, the HLL-Riemann solver finite volume method on unstructured triangular grids (Aliparast 2009, Yoon and Kang 2004, Harten, Lax and Van Leer 1983) is implemented for solving the resulting nonlinear-coupled system of SWEs and the Green-Ampt infiltration model. This infiltration model leads to nonlinear system, which is solved using implicit iterative methods. Appropriate discretization techniques are used for the bottom topography and Manning friction source terms to ensure the underlying physical properties. Moreover, piecewise linear reconstructions of the solutions are used at the discrete level in order to achieve the second-order accuracy of the numerical scheme. For subsurface water, we proposed accurate method for computing the infiltration rate under rainfall-runoff processes. The developed numerical techniques can be used to predict surface and subsurface water flows under rainfall events over complex bottom topography involving wet and dry areas and large bed slop where we obtain accurate results in terms of mass balance.

The paper is organized as follows. In Section 2, we present the coupled model of surface and subsurface water flows based on the SWEs and Green-Ampt model for one- and two-layers soil. To numerically solve this model, we propose to use the second-order HLL finite volume scheme in Section 3. In Section 4, we test the robustness and accuracy of the proposed numerical techniques for simulating rainfall-runoff over variable bottom topography with wet and dry areas where we perform three numerical experiments. Finally, some concluding remarks are given in Section 5.

## 2    MODEL EQUATION

In this study, we focus on the following shallow water system over variable bottom topography with different source terms such as friction effect, precipitation and infiltration, modelling coupled surface and subsurface water flows:

$$[1] \begin{cases} \dfrac{\partial h}{\partial t} + \dfrac{\partial hu}{\partial x} + \dfrac{\partial hv}{\partial y} = R - I, \\ \dfrac{\partial hu}{\partial t} + \dfrac{\partial}{\partial x}\left(hu^2 + \dfrac{g}{2}h^2\right) + \dfrac{\partial}{\partial y}(huv) = -gh\dfrac{\partial B}{\partial x} - \dfrac{\tau_x}{\rho}, \\ \dfrac{\partial hv}{\partial t} + \dfrac{\partial}{\partial x}(huv) + \dfrac{\partial}{\partial y}\left(hv^2 + \dfrac{g}{2}h^2\right) = -gh\dfrac{\partial B}{\partial y} - \dfrac{\tau_y}{\rho}, \end{cases}$$

where $h$ represents the water depth, $\mathbf{u} := (u,v)^T$ is the depth averaged velocity field of the flow, $g$ is the gravitational constant, the function $B(x,y)$ represents the bottom elevation, and $\rho$ is the water density. The source terms $R$ and $I$ are respectively the rainfall intensity and the infiltration rate into soil. We consider the following Manning formula for the computation of the friction source terms:

$$[2] \begin{cases} \dfrac{\tau_x}{\rho} = g\dfrac{n^2}{h^{1/3}}||\mathbf{u}||u, \\ \dfrac{\tau_y}{\rho} = g\dfrac{n^2}{h^{1/3}}||\mathbf{u}||v, \end{cases}$$

with $n$ is the Manning coefficient and $||\mathbf{u}||$ is the norm of the velocity field of the flow. The Manning formula is based on the assumption of uniform flow (Marcus, et al. 1992) and it is initially developed for flow in open channels and pipes. The Manning formula is widely used in many applications to compute the friction stresses in runoff models (Cea and Vazquez-Cendon 2012, Fernández-Pato, Gracia and García-Navarro 2018). In this study, we assume a constant Manning coefficient. While with this assumption we obtain a good agreement between experiments and simulations for the numerical tests performed in this study, a Manning formula with variable coefficient $n$ could be used, as shown in the comparative study by (Mugler, et al. 2011).

## 2.1 Green-Ampt infiltration model

In this section, we introduce the Green-Ampt infiltration model (Green and Ampt 1911) which is widely used (Tatard, et al. 2008, Mugler, et al. 2011) to compute the infiltration capacity of water in soils. This model is a simplified form of the Richards equation where the wetting front is considered as a sharp boundary between saturated and dry zones of the soil. To design an efficient coupled surface-subsurface numerical model in terms of computational cost, the Green-Ampt model is used in our study because of its simplicity and performance in computing the cumulative infiltration in soils.

We assume that water is ponded at the surface with depth $h_p$, the Green-Ampt model is used to compute the infiltration capacity $I_r$ as follows:

$$[3] \quad I_r = K_s \left( \frac{(\psi + h_p)\Delta\theta}{I_c} + 1 \right),$$

where $K_s$ is the saturated hydraulic conductivity of soil, $\Delta\theta = \theta_s - \theta_i$ with $\theta_s$ and $\theta_i$ are respectively the initial and saturated moisture contents, $\psi$ is the suction head, $I_c = \int_0^t I_r(\tau)\, d\tau$ is the cumulative infiltrated water. The wetting front reaches the depth $d_f = I_c/\Delta\theta$ at time $t$ where we assume $d_f = 0$ at $t = 0$. By integrating Eq. 3 over $[t, t + \Delta t]$ we obtain the following equation for the cumulative infiltration:

$$[4] \quad I_c^{t+\Delta t} = I_c^t + K_s \Delta t + \Delta\theta(\psi + h_p) \ln\left( \frac{(\psi + h_p)\Delta\theta + I_c^{t+\Delta t}}{(\psi + h_p)\Delta\theta + I_c^t} \right).$$

## 2.2 Two layers Green-Ampt infiltration model

Similarly to one-layer soil, we extend the Green-Ampt model for two-layers soil where the upper layer has a thickness $d_1$ with hydraulic conductivity $K_s = K_1$, suction head $\psi = \psi_1$ and $\Delta\theta = \Delta\theta_1$, the lower layer corresponds to the ground below the first one, with hydraulic conductivity $K_s = K_2$, suction head $\psi = \psi_2$ and $\Delta\theta = \Delta\theta_2$. The infiltration capacity for this system of two layers soil can be expressed as follows:

$$[5] \quad I_r = K_e \left( \frac{(\psi + h_p)}{d_f} + 1 \right),$$

where $K_e$ is the effective hydraulic conductivity given by:

$$[6] \quad K_e = \begin{cases} K_1, & \text{if} \quad d_f \leq d_1, \\ \dfrac{d_f}{\dfrac{d_1}{K_1} + \dfrac{d_f - d_1}{K_2}}, & \text{if} \quad d_f > d_1. \end{cases}$$

The cumulative infiltration for two layers soil denoted by $\tilde{I}_c$ is computed as:

$$[7] \quad \tilde{I}_c = I_c + a,$$

where the cumulative infiltration $I_c = d_f \Delta\theta$ is obtained using the following equation:

$$[8] \quad I_r = \frac{dI_c}{dt}.$$

By integrating over $[t, t + \Delta t]$, we get:

$$[9] \quad I_c^{t+\Delta t} = I_c^t + K_s \Delta t + A \times \ln\left( \frac{(\psi + h_p)\Delta\theta + I_c^{t+\Delta t}}{(\psi + h_p)\Delta\theta + I_c^t} \right),$$

where $K_s$, $\psi$ and $\Delta\theta$ are given according to the soil layer. The constants $a$ and $A$ are given by:

$$[10] \begin{cases} a = 0, \quad A = \Delta\theta_1(\psi + h_p), & \text{if} \quad d_f \leq d_1, \\ a = d_1\Delta\theta_1, \quad A = \Delta\theta_2\left(\psi + h_p - d_1\left(\frac{K_2}{K_1} - 1\right)\right), & \text{if} \quad d_f > d_1. \end{cases}$$

## 3 NUMERICAL METHOD

In this section, we derive the numerical method for the coupled model of surface and subsurface water flows. In order to satisfy the well-balanced property of the proposed numerical method, we used the vector of solution $U = (w, p, q)$ with new variables of the system. The variable $w = h + B$ is the water surface elevation, $p = hu$ and $q = hv$ are respectively the discharges in $x-$ and $y-$directions. By expressing the system in Eq. 1 with vector variable $U$, we obtain its following matrix form:

$$[11] \quad \frac{\partial U}{\partial t} + \frac{\partial H_x(U, B)}{\partial x} + \frac{\partial H_y(U, B)}{\partial y} = S_b(U, B) + S_f,$$

where the vector of fluxes $F = (H_x, H_y)^T$ and source terms $S_b$ and $S_f$ are:

$$[12] \begin{aligned} H_x &= \left(p, \frac{p^2}{w - B} + \frac{g}{2}(w - B)^2, \frac{pq}{w - B}\right)^T, \\ H_y &= \left(q, \frac{pq}{w - B}, \frac{q^2}{w - B} + \frac{g}{2}(w - B)^2\right)^T, \\ S_b &= \left(I - R, -g(w - B)\frac{\partial B}{\partial x}, -g(w - B)\frac{\partial B}{\partial y}, 0\right)^T, \\ S_f &= \left(0, -\frac{\tau_x}{\rho}, -\frac{\tau_y}{\rho}, 0\right)^T. \end{aligned}$$

### 3.1 Surface water computation

We assume that the domain is partitioned into a set of triangular cells $E_j$ of area $|E_j|$. Let $E_{jk}$ $k = 1,2,3$ be the neighbouring cells of $E_j$ with common edges $e_{jk}$ of length $l_{jk}$, and define $\mathbf{n}_{jk} = (\cos(\theta_{jk}), \sin(\theta_{jk}))$ the unit vector normal to the edge $e_{jk}$ and points toward the cell $E_j$. Denote by $G_j = (x_j, y_j)$ the coordinates of the center of mass of the cells $E_j$ and $P_{jk} = (x_{jk}, y_{jk})$ the midpoint of the edge $e_{jk}$ with corresponding vertices $P_{jki}$ $i = 1,2,3$, as shown in Figure 1.

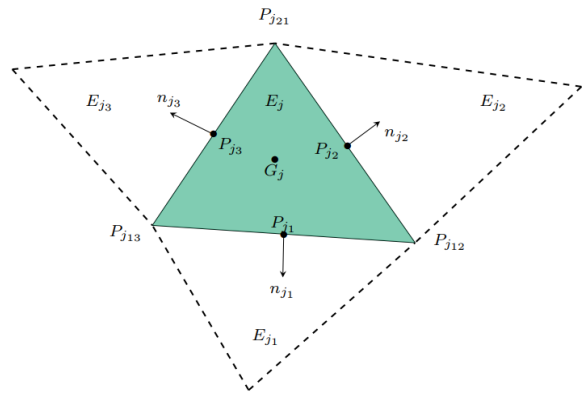

Figure 1: The triangular grids.

The second-order HLL-Riemann solver schemes on unstructured triangular grid (Yoon and Kang 2004) is applied to the coupled model of surface and subsurface water flows, which can be written as follows:

[13] $\dfrac{\partial \overline{U}_j}{\partial t} = \dfrac{1}{|E_j|} \sum_{k=1}^{3} F_{jk} \cdot n_{jk} l_{jk} + \overline{S}_{b_j} + \overline{S}_{f_j},$

where $\overline{U}_j = \dfrac{1}{|E_j|} \int_{E_j} U(t,x,y) dx dy$ is the approximation of the cell average of the computed solution, and $\overline{S}_{b_j} = \dfrac{1}{|E_j|} \int_{E_j} S_b(U,B) dx dy$ and $\overline{S}_{f_j} = \dfrac{1}{|E_j|} \int_{E_j} S_f(t,x,y) dx dy$ are the cell averages of the source terms.

The advective numerical fluxes $F_{jk}$ are:

[14] $F_{jk} \cdot n_{jk} = \begin{cases} F_{jk}^L(U_{jk}^L, B_{jk}) \cdot n_{jk}, & \text{if } S_{jk}^L \geq 0, \\ \dfrac{S_{jk}^R F_{jk}^L(U_{jk}^L, B_{jk}) \cdot n_{jk} - S_{jk}^L F_{jk}^R(U_{jk}^R, B_{jk}) \cdot n_{jk} + S_{jk}^L S_{jk}^R (U_{jk}^R - U_{jk}^L)}{S_{jk}^R - S_{jk}^L}, & \text{if } S_{jk}^L \leq 0 \leq S_{jk}^R, \\ F_{jk}^R(U_{jk}^R, B_{jk}) \cdot n_{jk}, & \text{if } S_{jk}^R \leq 0. \end{cases}$

where $S_{jk}^L$ and $S_{jk}^R$ are respectively the left and right values of wave speeds of propagation, which are estimated in our study using the extreme eigenvalues based on the left and right values of the water depth and velocities. The values $U_{jk}^L := U_j(P_{jk})$ and $U_{jk}^R := U_{jk}(P_{jk})$, correspond respectively to the left and right reconstructed values of the vector variable $U$ at the midpoints $P_{jk}$, which are obtained using the following linear approximation:

[15] $\widetilde{U}(x,y) = \overline{U}_j + (U_x)_j (x - x_j) + (U_y)_j (y - y_j).$

In order to ensure the positivity of the water depth $h$, this reconstruction is corrected for the water surface elevation $w$ by introducing the parameter $\alpha \in [0,1]$ (Beljadid, Mohammadian and Kurganov 2016):

[16] $\widetilde{w}(x,y) = \overline{w}_j + \alpha(w_x)_j (x - x_j) + \alpha(w_y)_j (y - y_j).$

The parameter $\alpha$ is computed by requiring $w_j(P_{jk}) > 0$ and $w_{jk}(P_{jk}) > 0$. The numerical gradients in the piecewise linear reconstruction Eqs. 15 and 16 are obtained using the Green-Gauss theorem:

[17] $\nabla U_j = \dfrac{1}{|E_j|} \int_{E_j} \nabla U_j(x,y) dx dy = \dfrac{1}{|E_j|} \sum_{s=1}^{3} \int_{E_j} \widetilde{U}_{js} n_{js} de,$

where the values $\widetilde{U}_{js}$ are estimated at the cell interface adopting the approach procedure developed in (Beljadid, Mohammadian and Kurganov 2016). The required values of the bottom topography at the midpoints $B_{jk}$ are computed from its known values at the vertices, obtained using a continuous piecewise linear approximation:

[18] $B_{jk} = \dfrac{B_{jk1} + B_{jk2}}{2}.$

We used linear reconstruction for the bottom topography at each computational cell. The estimated values of the bottom topography at the midpoints are used to compute the elevation at the center of mass of the cell:

[19] $B_j = \dfrac{B_{j1} + B_{j2} + B_{j3}}{3},$

## 3.2 Discretization of the source terms

In order to guarantee the well-balance property of the numerical scheme, appropriate discretization of the cell averages of the source term $\overline{S}_{b_j}$ is needed, to exactly balance the numerical fluxes. To this end, the well-balance discretization of the source term $\overline{S}_{b_j}$ is given by (Bryson, Epshteyn and Kurganov 2011):

$$[20] \quad \begin{aligned} \bar{S}_{b_j}^{(1)} &= R_j - I_j, \\ \bar{S}_{b_j}^{(2)} &= \frac{g}{2|E_j|} \sum_{k=1}^{3} l_{jk}\left(w_{jk}^L - B_{jk}\right)^2 \cos(\theta_{jk}) - g(w_x)_j (\overline{w}_j - B_{jk}), \\ \bar{S}_{b_j}^{(3)} &= \frac{g}{2|E_j|} \sum_{k=1}^{3} l_{jk}\left(w_{jk}^L - B_{jk}\right)^2 \sin(\theta_{jk}) - g(w_y)_j (\overline{w}_j - B_{jk}). \end{aligned}$$

For the discretization of the Manning's friction term $\overline{S}_{f_j}$, we used the methodology developed in (Xia and Liang 2018).

## 3.3 Sub-surface water computation

In the previous section we used the Green-Ampt equation describing the infiltration of water flow into soil where the water is available at the surface at all times. However, during rainfall event the water ponded at the surface if the rainfall intensity is greater than the infiltration capacity of the soil, otherwise all the water infiltrates into the soil. In this study, the infiltration rate is computed using the following techniques:

- If $R_j^n > I_{r_j}^n$ : water is ponded at the surface. The infiltration rate is computed from Eqs. 3 to 10 by:

$$[21] \quad I_j^n = \min\left(\frac{\bar{h}_j^n}{\Delta t}, I_{r_j}^n\right).$$

- If $R_j^n \leq I_{r_j}^n$ : no ponding occurs at the surface, then we compute the infiltration rate by:

$$[22] \quad \begin{cases} I_j^n = R_j^n, \\ I_{c_j}^{n+1} = I_{c_j}^n + I_j^n \times \Delta t. \end{cases}$$

## 4 NUMERICAL EXPERIMENTS

In this section, we will test the performance of the proposed numerical scheme and its ability to simulate surface and subsurface water flows. In our simulations, we used the gravitational constant $g = 9.81$. In the first example, we perform a validation of our numerical techniques using experimental data for a simple laboratory non infiltrating slop surface. The second example focuses on the well-balanced and conservative properties of the scheme where we perform a numerical test over infiltrating surface using one-and two-layers soil. In the last example, we test the capability of the proposed numerical method to simulate surface-subsurface water flow over complex bottom topography.

### 4.1 Rainfall-runoff over a slop surface

This numerical example is considered to validate the proposed numerical techniques using experimental data (Gottardi and Venutelli 2008). The computational domain is a plane of length $L = 21.945$ m and width $l = 1$ m with a slop $S_x = 0.04$ m which is discretized using triangular grids with an average cell area $|\bar{E}_j| = 3.9 \times 10^{-3}$. The soil surface is initially dry and the rainfall with an uniform intensity $R$ is imposed over the domain. The Manning's coefficient $n = 0.5$ is used by (Gottardi and Venutelli 2008). Hydrographs of computed and measured outlet discharges are presented in Figure 2 (left) for different values of the Manning coefficient to show the impact of uncertainty of this coefficient on the results. A suitable value of the parameter $n$ should be used to calibrate the proposed numerical model. We first calibrate the model where we compare the computed and measured outlet discharges to obtain the appropriate value of the Manning coefficient by minimizing the root mean square error (RMSE) defined as follows:

[23] $RMSE = \sqrt{\frac{1}{N}\sum_{k=1}^{N}(q_k^{(s)} - q_k^{(o)})^2}$,

where $q_k^{(s)}$ and $q_k^{(o)}$ are the simulated and observed outlet discharges, and $N$ is the total number of observations. We obtain $n = 0.48$ for the calibrated model and Figure 2 (left) shows the corresponding results. To validate the calibrated model against other experimental data, numerical simulations are performed using another value of rainfall intensity as shown in Figure 2 (right). We observe good overall agreement between experiments and simulations.

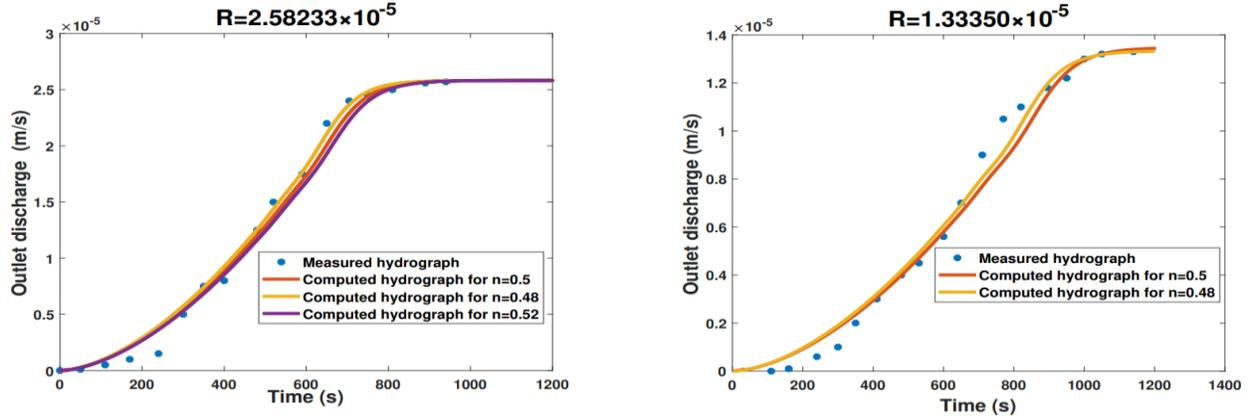

Figure 2: Hydrographs of computed and measured outlet discharges for constant rainfalls over slope surface.

## 4.2 Well-balanced and conservative properties

In this example, we test the well-balanced property and the conservation of the proposed numerical scheme. The proposed method is applied to solve the system Eq.1 for one-and two-layers soil, where we considered the Silt loam soil of 30 percent initial effective saturation and Sandy loam soil of initial effective saturation of 40 percent. The Green-Ampt parameters of the used soils are given in Table 1 (Te Chow, Maidment and Mays 1988). The computational domain $[0,2] \times [0,1]$ is discretized using triangular grids with an average cell area $|\bar{E}_j| = 6.3735 \times 10^{-4}$. The surface water is initially constant everywhere with depth $h(0) = 0.1$ and zeros field velocity $\mathbf{u} = 0$. The wall conditions are applied at the boundaries of the domain.

Table 1: Green-Ampt parameters.

|  | Soil | Suction at the wetting front $\psi$ (cm) | Saturated hydraulic conductivity $K_s$ (cm/h) | $\Delta\theta = \theta_s - \theta_i$ | $d_1$ (cm) |
|---|---|---|---|---|---|
| Upper layer | Sandy loam | 11.01 | 1.09 | 0.247 | 0.1 |
| Lower layer | Silt loam | 16.7 | 0.65 | 0.340 |  |

Figure 3, represents the cross section along the $x$-axis of the computed water depth at different times over infiltrating surface and the corresponding total mass balance for one-and two-layers soil. The results show that the well-balanced property of the scheme is satisfied where water surface elevation decreases with time and remains constant in space. Moreover, the computed total mass of water which includes the surface water and the cumulative infiltrated water is preserved with time, that is $h(t) + I_c(t) = 0.1 = h(0)$.

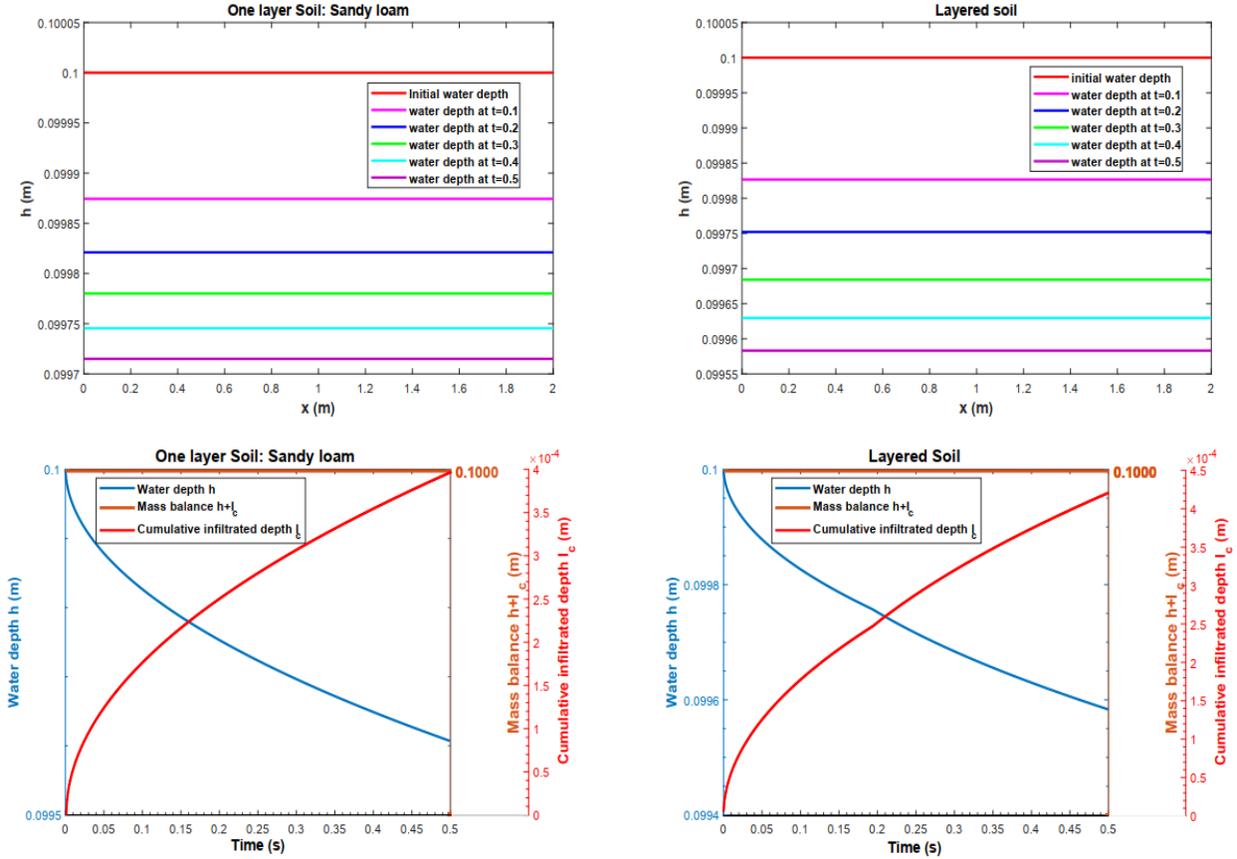

Figure 3: Computed water depth and the corresponding total mass for one layer soil: Sandy loam (left) and two-layers soil (right).

### 4.3 Surface-subsurface flows over complex basin

We consider a 2D modified numerical example used in (Park, Kim and Kim 2019), to test the ability of the proposed scheme for simulating runoff processes. In this example, we used the computational domain $[0,10] \times [0,8]$ and the following variable topography:

$$[23]\quad B(x,y) = 0.01y + 0.01|x - 0.5| - 0.01 \sin\left(\frac{\pi x}{2}\right) - 0.01 \sin\left(\frac{\pi y}{2}\right) + 1.$$

The Manning's roughness coefficient is set to $n = 0.013\ m^{-1/3}s$, and wall boundary conditions are imposed at the sides and upstream boundaries, while outflow condition is used at the downstream of the domain. We used the hydraulic conductivity of the soil $K_s = 7$ mm/h, suction head $\psi = 50$ mm and $\Delta\theta = 0.125$. A constant rainfall with intensity $R = 500$ mm/h is applied over the domain for 5 min.

Figure 4 illustrates the evolution of the water surface elevation $w$, the water depth $h$ and the component of the velocity in y-direction at different times. Initially, the domain is dry and at time $t = 1$ min, the water starts ponding at the surface and filled ponds in the upper areas, the excessive water overflows the ponds and the lower areas inside the domain. At time $t = 5$ min, the rain stops while water keeps moving at the surface to the downstream of the domain, and most areas became flooded with time. At time $t = 8$ min the upper areas of the domain become almost dry except inside the ponds. In Figure 5, we show the flow field velocity vectors of the computed solution at time $t = 5$ min. We observe a rapid increase of the flow velocity in the downstream. Finally, we conclude that our numerical simulations using the proposed numerical techniques reproduce physically reasonable rainfall-runoff processes over complex bottom topography with wet and dry area and large bed slop.

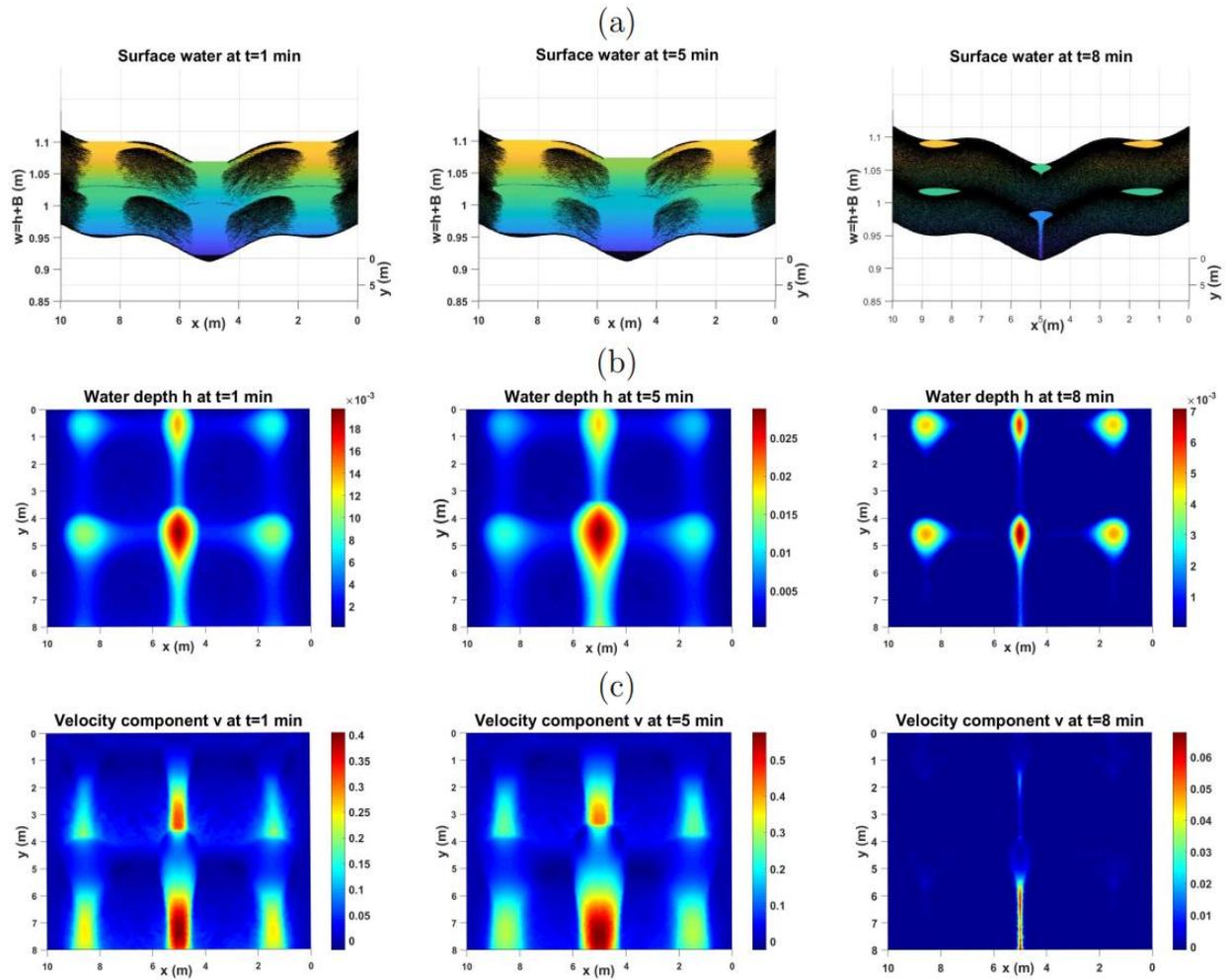

Figure 4: Evolution of the computed surface water (a), water depth (b) and component of the velocity in y-direction (c) at times $t = 1$ min, $t = 5$ min and $t = 8$ min.

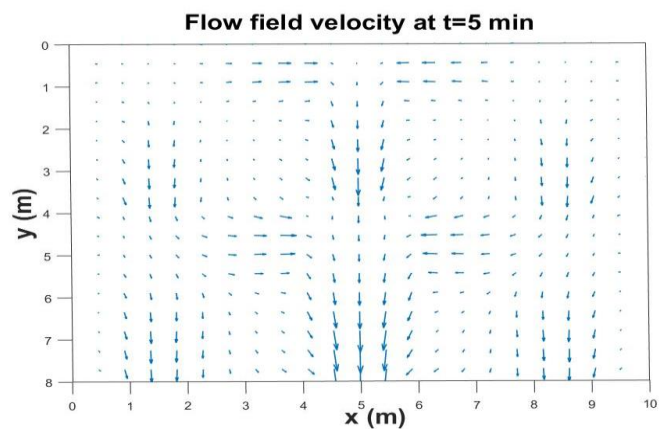

Figure 5: Flow field velocity vector of the computed solution at time $t = 5$ min.

## 5 CONCLUSION

We proposed numerical techniques for modelling coupled surface and subsurface water flows. The shallow water equations are used for surface water flow over variable bottom topography with source terms due to friction effect, precipitation and infiltration. For subsurface water flow, we used the Green-Ampt infiltration model with one-and two-layers soil. HLL Riemann solver finite volume schemes with linear reconstruction of the solutions are implemented for solving the coupled nonlinear system. Appropriate discretization techniques are used for the bottom topography and Manning friction source terms to ensure the underlying physical properties of the proposed numerical scheme. Moreover, we used accurate techniques for subsurface water computation under rainfall-runoff processes. Our numerical results confirm the accuracy and the capability of the proposed numerical techniques for modelling surface-subsurface water flows over variable bottom topography involving wet and dry areas and large bed slop.